\documentclass[11pt]{article}

\usepackage[T1,T2A]{fontenc}
\usepackage[utf8]{inputenc}
\usepackage{amsfonts}
\usepackage{amssymb}
\usepackage{epsfig}
\usepackage{tikz}
\usepackage{comment}
\usepackage{epstopdf}%added by Gareth

\usepackage{theorem}

\usepackage{color} % added by Gareth

\newtheorem{theorem}{Theorem}[section]

{\theorembodyfont{\rmfamily} % No need in ``\rm'' command after ``\begin''

\newtheorem{remark}[theorem]{Remark}
\newtheorem{example}[theorem]{Example}

}

\def\eq{\,=\,}
\def\Z{\mathcal{Z}} % added by Gareth
\def\C{\mathbb{C}}  % added by Gareth
\def\Q{\mathbb{Q}}
\def\F{\mathbb{F}}
\def\R{\mathcal{R}} % changed by Gareth
\def\D{\mathcal{D}} % added by Gareth
\def\X{\mathcal{X}} % added by Gareth
\def\Y{\mathcal{Y}} % added by Gareth
\def\P{\mathbb{P}}  % added by Gareth

\def\msk{\medskip}
\def\ssk{\smallskip}

\def\D{\mathcal{D}} % added by Gareth

\title{Hurwitz groups as monodromy groups of dessins: several examples}
\date{\today}
\author{Gareth A. Jones and Alexander K. Zvonkin}

\begin{document}

\maketitle

\begin{abstract}
We present a number of examples to illustrate the use of small quotient dessins as substitutes for their often much larger and more complicated Galois (minimal regular) covers. In doing so we employ several useful group-theoretic techniques, such as the Frobenius character  formula for counting triples in a finite group, pointing out some common traps and misconceptions associated with them. Although our examples are all chosen from Hurwitz curves and groups, they are relevant to dessins of any type.
\end{abstract}

\noindent{\bf MSC Classification}: primary 14H57, secondary 20B25.

\noindent{\bf Key words}: Hurwitz curve, Hurwitz group, dessin d'enfant, automorphism group, monodromy group.

%%%%%%%%%%%%%%%%

\section{Introduction}

When Grothendieck wrote his Esquisse d'un Programme, he famously expressed 
his delight that basic objects such as plane trees, and other simple drawings, 
could encode highly sophisticated mathematical structures, namely algebraic curves 
defined over number fields. Indeed, such was his pleasure in the elegance and 
power of this data compression that he used the childish-looking nature of these 
sketches in giving them the name of {\em dessins d'enfants}, which has 
now come to represent a vast mathematical theory (see~\cite{GG-D-12, JW-16, LZ-04} for instance). 

The recent proof by Gonz\'alez-Diez and Jaikin-Zapirain~\cite{GDJZ-15} that 
regular dessins provide a faithful representation of the absolute Galois group 
has focussed attention even more closely than before on these highly symmetric 
objects, those dessins for which the associated Bely\u\i\/ function is a regular 
covering. Every dessin $\D$ is the quotient of some regular dessin 
$\R=\widetilde\D$ where $\widetilde\D$ is the Galois (minimal regular) cover of $\D$, 
by a subgroup $H$ of the automorphism group $G$ of $\R$; this latter group is 
also realised as the monodromy group of $\D$, or equivalently of its associated 
Bely\u\i\/ function, regarded as a branched covering of the sphere; this is a 
permutation group acting transitively on the fibre over a base-point, with 
$H$ as a point-stabiliser.

%%%%%%%%%%%%%%%%%%%%%%%%%%%%%%%%%%%%%%%%%%%
\def\softd{{\leavevmode\setbox1=\hbox{d}%
\hbox to 1.05\wd1{d\kern-0.4ex{\char039}\hss}}}%cstocs
\def\softt{{\leavevmode\setbox1=\hbox{t}%
\hbox to \wd1{t\kern-0.6ex{\char039}\hss}}}%cstocs
\def\softl{l\kern-0.45ex\raise0.1ex\hbox{'}\kern-0.10ex}%cstocs
\def\softL{L\kern-0.8ex\raise0.1ex\hbox{'}\kern0.1ex}%cstocs
%%%%%%%%%%%%%%%%%%%%%%%%%%%%%%%%%%%%%%%%%%%

A common situation is that in which a dessin $\D$, which is `small' in some sense 
(having low genus, or few edges, for example) is used in this way as a substitute 
for a much larger regular dessin $\R$. The latter is uniquely determined as the 
Galois cover of $\D$ provided $G$ acts faithfully on the cosets of $H$, that is, 
the core of $H$ in $G$ is trivial, in which case we will call $\D$ a {\em faithful\/} 
quotient of $\R$. For instance, one can see this idea in action, even before 
the era of dessins, in the use by Conder~\cite{Conder-80, Conder-81} of Graham 
Higman's technique of `sewing together coset diagrams' to realise large alternating 
and symmetric groups as quotients of triangle groups. Similarly, Jendro\softl, 
Nedela and \v Skoviera~\cite{JeNeSk-97} have used this idea  to obtain new results 
and new proofs of old ones for graphs and for maps on surfaces.

Our aim in this paper is to explore this relationship between regular dessins $\R$ 
and their faithful quotients $\D$, with its occasional unexpected subtleties. 
We do this through a series of examples, chosen with the secondary aim of 
illustrating some useful techniques from group theory (finite, discrete and 
computational) for enumerating, constructing and classifying dessins of various 
types associated with specific groups. Our examples are all Hurwitz groups and 
curves, those attaining Hurwitz's upper bound~\cite{Hurwitz-1893} of $84(g-1)$ 
for the number of automorphisms of a curve of genus $g \ge 2$. This choice 
is purely for personal and historical interest 
since the corresponding dessins are extremely rare. (Conder~\cite{Conder-homepage} listed all the regular dessins of genus 2 to 101; their total number is 19\,029, and only seven of them attain the Hurwitz bound.) Nevertheless, the methods we describe can in fact be applied to dessins of any type.

There are no new theorems in this paper. Indeed, the only proof we offer is 
really a disproof, of the occasional and erroneous assertion that the Hurwitz 
group of genus $17$ is isomorphic to the affine group ${\rm AGL}_3(2)$. 
These two groups do indeed look very similar, both being extensions of an 
elementary abelian normal subgroup of order $8$ by ${\rm GL}_3(2)$; however, 
the affine group is a split extension, while the Hurwitz group is not.
We also exhibit certain traps to avoid while using irreducible characters and the Frobenius formula to count the number of dessins with a given monodromy group.
     
We will use many well-known facts concerning the inner structure of various
finite groups. Unfortunately (and obviously), we are unable to fill in all the 
necessary details, since otherwise the paper would be enlarged {\em ad infinitum}. 
In a few occasions we supply the reader with some details but in most cases he
or she should consult other sources on group theory.

\msk

\noindent{\bf Acknowledgement} The authors gratefully thank Jarke van Wijk for giving permission to use his beautiful pictures of maps (Figures~\ref{fig:Klein} and \ref{fig:Fricke-Macbeath}), and the referee for some very helpful suggestions. The second author is partially supported by the research grant {\sc Graal} ANR-14-CE25-0014.

%%%%%%%%%%%%%%%%%%%%%%%%%%%%

\section{Background details}

\subsection{Dessins}

By Bely\u\i's Theorem, a compact Riemann surface $\mathcal S$, regarded as 
a complex projective algebraic curve, is defined over a number field if and 
only if it admits a Bely\u\i\/ function, a non-constant meromorphic function 
$f:{\mathcal S}\to\P^1(\C)$ ramified over at most three points (which one can, 
without loss of generality, take to be $0, 1$ and $\infty$). The 
{\em monodromy group\/} $G$ of the Bely\u\i\/ pair $(\mathcal S,f)$ is the 
monodromy group of the covering of the thrice-punctured sphere 
$\P^1(\C)\setminus\{0, 1, \infty\}$ induced by $f$, that is, the permutation 
group on the sheets (more precisely on the fibre over a base-point) induced by 
unique lifting of closed paths. If $f$ has degree $n$ then $G$ is a transitive 
subgroup of the symmetric group ${\rm S}_n$, generated by the local monodromy 
permutations $x, y$ and $z$ around $0, 1$ and $\infty$, satisfying $xyz=1$. 
By the Riemann Existence Theorem, these three permutations define $\mathcal S$ 
and $f$ up to isomorphism.

Following Grothendieck, one can represent $f$ by means of a bicoloured map 
on $\mathcal S$, called a dessin $\D$, with the fibres over $0$ and $1$ as the 
black and white vertices,
and the unit interval lifting to the $n$ edges, one on each sheet of the covering. 
Then the permutations $x$ and $y$ represent the rotations of the edges around their 
incident black and white vertices, while the permutation $z=(xy)^{-1}$ rotates edges, 
two steps at a time, around incident faces.

The {\em automorphism group\/} ${\rm Aut}(\D)$ of a dessin $\D$ is the group 
of covering transformations of $f$, or equivalently the centraliser of $G$ in 
${\rm S}_n$. A dessin $\D$ is {\em regular\/} if the covering is regular, 
that is, if ${\rm Aut}(\D)$ acts transitively on the edges, in which case 
${\rm Aut}(\D)\cong G$.

The {\em type} of a dessin, or of a triple $(x, y, z)$, is the triple of orders 
of $x, y$ and $z$. Rather more information is conveyed by the {\em passport\/}, 
by which we mean the ordered triple of partitions of $n$ giving the cycle-structures 
of $x, y$ and $z$. These correspond to the conjugacy classes of ${\rm S}_n$ containing 
these permutations, but in some cases we will give more precise information by 
referring to the passport as the triple of their conjugacy classes in the monodromy 
group $G$ which they generate.

In cases where $y^2=1$ the white vertices of the dessin are redundant and may be 
omitted, so that fixed points of $y$ are now represented by free edges 
(one end of which is a vertex while the other remains free),
and 2-cycles by traditional edges or loops. The edges of the original dessin now 
correspond to the half-edges of the resulting map, and the monodromy group of the dessin can be identified with that of the map, now permuting half-edges. In some of our examples, 
we will use this simplification without further comment.

For further background reading we suggest~\cite{GG-D-12, JW-16} or~\cite{LZ-04}.

%%%%%%%%%%

\subsection{The Frobenius formula}

For a given finite group $G$, the regular dessins of type $(p,q,r)$ with automorphism 
group $G$ are in bijective correspondence with the torsion-free normal subgroups 
$N$ of the triangle group
\[\Delta=\Delta(p,q,r)=\langle  X, Y, Z\mid X^p=Y^q=Z^r=XYZ=1\rangle\]
with $\Delta/N\cong G$. These subgroups $N$ correspond bijectively 
to the orbits of ${\rm Aut}(G)$ on generating triples $(x, y, z)$ of elements 
of orders $p, q$ and $r$ in $G$ satisfying $xyz = 1$. Since the action of ${\rm Aut}(G)$
on generating triples is semiregular, meaning that only the identity element has fixed points, the number of orbits, and hence of dessins, is equal to the number 
of such triples divided by $|{\rm Aut}(G)|$. The first step in the calculation 
of this number is to use the following classical result~\cite{Frobenius-1896}
(for a modern treatment see \cite[Ch.~7]{Serre}):

\begin{theorem}[Frobenius]\label{th:frobenius}
Let $\X$, $\Y$ and $\Z$ be conjugacy classes in a finite group $G$. Then the number of solutions in $G$ of the equation $xyz=1$, where $x\in\X$, $y\in\Y$ and $z\in\Z$, is given by the formula
\[\frac{|\X|\cdot |\Y|\cdot |\Z|}{|G|}
\sum_{\chi}\frac{\chi(x)\chi(y)\chi(z)}{\chi(1)}, \]
where the sum is over all irreducible complex characters $\chi$ of $G$.
\end{theorem}

\begin{remark}
The first part of this formula, omitting the character sum, can be regarded as a naive guess for the number of solutions, assuming that the values of $xyz$ are evenly distributed over the elements of $G$. Of course, in general they are not, and the character sum can be regarded as a correction term, taking into account the particular structure of the chosen group $G$. In many cases the character sum is dominated by the contribution, equal to $1$, from the principal character, in which case the naive guess is not far from the correct answer.
\end{remark}

\begin{remark}
This is a particular case of a more general formula for the number of solutions of $x_1\ldots x_k=1$ in $G$, where each $x_i$ is chosen from a specific conjugacy class. The only part of the generalisation which is not obvious is the denominator in the character sum, which is $\chi(1)^{k-2}$. For this, and other similarly useful formulae, see~\cite[Ch.~7]{Serre}.
\end{remark}

\begin{remark}
The {\sc Atlas} of Finite Groups~\cite{Atlas-05} contains character tables of many 
finite simple groups, and of other associated groups. Maple can compute characters 
of symmetric groups, and GAP can compute character tables of arbitrary not-so-big
groups.

A conjugacy class $\X$ in $G$ has order
\[|\X]=\frac{|G|}{|C_G(x)|}\]
where $C_G(x)$ is the centraliser in $G$ of an element $x\in\X$, and similarly, 
for $\Y$ and $\mathcal Z$. The {\sc Atlas} gives orders of centralisers, 
rather than those of conjugacy classes. 

\end{remark}

The Frobenius formula gives us the number of triples in $G$ with passport $(\X, \Y, \mathcal Z)$. To obtain the number of triples of type $(p,q,r)$, one simply takes the sum of these numbers over all triples of conjugacy classes $\X, \Y$ and $\mathcal Z$ consisting of elements of orders $p, q$ and  $r$.

Suppose we are given a triple $\pi=(\lambda,\mu,\nu)$ of partitions of number $n$, 
and we would like to know the number of dessins with the passport $\pi$. Enumerative combinatorics gives an explicit answer only in very specific cases, 
such as, for example, for plane trees. In more complicated cases, an invaluable, and 
in most cases the only source of information is the Frobenius formula. But it must be 
used with care.

Let us discuss first the case of the symmetric group $G={\rm S}_n$. Since in this 
case a cycle-structure uniquely determines the corresponding conjugacy class, 
it seems that we can then apply the Frobenius formula directly. There is, however, 
a trap to avoid: a triple of permutations $(x,y,z)$ with cycle-structures 
$(\lambda,\mu,\nu)$ does not necessarily generate a transitive subgroup of 
${\rm S}_n$. We must find a way to eliminate these non-transitive solutions.

The next difficulty to resolve is the fact that the edges of the same dessin~$\D$ 
may be labelled in many different ways. To be specific, the number of labellings 
is $n!/|{\rm Aut}(\D)|$. Therefore it is reasonable to divide the number of triples
of permutations by $n!$: in this way we will get the ``number'' of non-isomorphic
dessins, each one of them being counted with the weight $1/|{\rm Aut}(\D)|$.

\begin{example}\label{ex:6321}
Let us take $n=12$ and look for those maps with the passport 
$\pi=(6^1 3^1 2^1 1^1, 2^6, 6^1 3^1 2^1 1^1)$.
Computing the corresponding characters using Maple, applying the Frobenius formula
and dividing the result by $12!$ we get $19\frac{1}{2}$. In fact, the correct
answer is 18. Two more ``maps'' are non-connected, and one of them has a
non-trivial automorphism of order 2, so that its contribution to the sum is $1/2$.
\end{example}

Now suppose that we work inside a group $G$ different from ${\rm S}_n$. Then, the
next difficulty arises: we can come across a cycle-structure corresponding to several
different conjugacy classes. For example, in the group ${\rm PSL}_2(27)$
(in its natural representation of degree $n=28$), which will be treated in 
Section~\ref{sec:genus-118}, there are two conjugacy classes $\X_1$, $\X_2$ of 
elements with cycle-structure $3^9 1^1$, one conjugacy class $\Y$ of elements with 
cycle-structure $2^{14}$, and three conjugacy classes $\Z_1$, $\Z_2$, $\Z_3$ with 
cycle-structure $7^4$. Thus, looking for dessins with passport 
$\pi=(3^9 1^1, 2^{14}, 7^4)$ we must take into account six possible combinations 
of conjugacy classes.

Return now to the problem with many possible labellings of edges of the same dessin.
If we want to stay inside a given group $G$ and not to be sent to one of its
conjugate copies, division by $n!$ would be not a good idea. We are tempted to
divide the number of triples $(x,y,z)$ by $|G|$. But there is a trap which 
awaits us here, and it is much more subtle than the previous ones! A significant 
part of our paper is devoted to untangling the complications arising in this 
case. See in this respect the discussions in \S\ref{sec:quotients} and
Sections~\ref{sec:genus-14} and \ref{sec:genus-118}.

Finally, a triple $(x,y,z)$ of permutations belonging to a group $G$ may generate
not the entire group $G$ but only a proper subgroup. If we are interested 
in triples generating $G$ itself we need to do more work depending on the particular
group and its structure. One possible technique for doing this is discussed in the next section.

%%%%%%%

\subsection{M\"obius inversion in groups}

Instead of counting all triples of a given type $(p,q,r)$ in $G$, we need to count {\em generating\/} triples of that type. Occasionally it is obvious which triples generate $G$ (all of them  in some cases), but in general a more systematic method, based on P.\,Hall's theory of M\"obius inversion in groups~\cite{Hall-36}, is available.

Given a type $(p,q,r)$, let $\sigma(H)$ denote the number of triples of that type in each subgroup $H\le G$, and let $\phi(H)$ denote the number of them which generate $H$. Since each triple in $G$ generates a unique subgroup, we have
\begin{equation}\label{sigmaeqn}
\sigma(G)=\sum_{H\le G}\phi(H).
\end{equation}
As shown by Hall, equation (\ref{sigmaeqn}) can be inverted, reversing the roles of the functions $\sigma$ and $\phi$, to give
\begin{equation}\label{phieqn}
\phi(G)=\sum_{H\le G}\mu_G(H)\sigma(H),
\end{equation}
where $\mu_G$ is the M\"obius function on the lattice of subgroups of $G$, recursively defined by
\begin{equation}\label{mueqn}
\sum_{K\ge H}\mu_G(K)=\delta_{H,G}\,,
\end{equation}
with $\delta$ denoting the Kronecker delta.

Although for many finite groups $G$, such as all but a few of the alternating and symmetric groups, the subgroup lattice is too complicated to allow $\mu_G$ to be calculated, this has been achieved for several classes of groups. For instance Hall dealt with nilpotent groups and the groups ${\rm PSL}_2(p)$ ($p$ prime) in~\cite{Hall-36}, and the latter calculation was extended by Downs~\cite{Downs-91} to ${\rm PSL}_2(q)$ and ${\rm PGL}_2(q)$ for all prime powers $q$. He and the first author have recently dealt with the Suzuki groups in~\cite{DJ-16}, and Pierro with the `small' Ree groups in~\cite{Pierro-14}. In many cases (for instance, if $H$ is not an intersection of maximal subgroups of $G$) we find that $\mu_G(H)=0$, so such subgroups $H$ can be omitted from the summation in equation~(\ref{phieqn}). By using the Frobenius formula to evaluate $\sigma(H)$ for subgroups $H\le G$ with $\mu_G(H)\ne 0$ one can determine $\phi(G)$, and hence obtain the number $\phi(G)/|{\rm Aut}(G)|$ of regular dessins of a given type with automorphism group~$G$.

The examples we present later in this paper do not in fact require M\"obius inversion, since in most cases it is easy to see that the relevant triples generate the whole group. However, there are more complicated cases, such as those considered in~\cite{DJ-16, Pierro-14}, where it cannot be avoided.

Hall's theory has much wider applications than that of counting triples described here. In its most general form, it can be used to count the normal subgroups of {\em any} finitely generated group with a given finite quotient group. For instance, Hall showed that the free group of rank 2 has 19 normal subgroups $N$ with quotient isomorphic to the alternating group ${\rm A}_5$; it follows that there are 19 regular dessins $\R$ with ${\rm Aut}(\R)\cong {\rm A}_5$ (described in~\cite{BJ-01}), and since this group has eight faithful transitive permutation representations, there are $19\cdot 8 = 152$ dessins $\D$ with monodromy group ${\rm A}_5$.

%%%%%%%

\subsection{Hurwitz groups and surfaces}

We have chosen all our examples from the historically important Hurwitz groups. 
A {\em Hurwitz group\/} $G$ is a finite group which attains Hurwitz's upper 
bound~\cite{Hurwitz-1893} of $84(g-1)$ for the order of the automorphism group 
of a compact Riemann surface of genus $g \ge 2$. Equivalently, $G$ is 
a non-trivial finite quotient $\Delta/N$ of the triangle group 
$\Delta=\Delta(2, 3, 7)$, acting as the automorphism group of the Riemann 
surface ${\mathcal S}={\mathbb H}/N$ (called a {\em Hurwitz surface} or 
{\em Hurwitz curve}), where $\mathbb H$ is the hyperbolic plane. 
As such, $G$ is the automorphism group of a regular dessin $\R$ on $\mathcal S$, 
called a {\em Hurwitz dessin}; regarded as a map, this is a trivalent tessellation 
by heptagons, or its dual, a 7-valent triangulation. The number of such surfaces 
and dessins associated with $G$, up to isomorphism, is equal to the number of normal 
subgroups $N$ of $\Delta$ with $\Delta/N\cong G$. Note that since its three periods 
are mutually coprime, $\Delta$ is a perfect group (that is, it has no non-trivial abelian quotient groups), and hence so is every Hurwitz group $G$. In particular, this implies that a Hurwitz group cannot be solvable.

The choice of the parameters $(p,q,r)=(2,3,7)$ is explained as follows. A regular 
dessin of type $(p,q,r)$ with an automorphism group $G$ of order $|G|=n$ has
$n/p$ black vertices of valency $p$, $n/q$ white vertices of valency $q$,
$n/r$ faces of valency $r$, and $n$ edges. Its Euler characteristic is thus equal to
$$
2-2g \eq \frac{n}{p}+\frac{n}{q}+\frac{n}{r}-n.
$$
In order to get a negative Euler characteristic with the least possible genus
(when $n$ is given) we should have the sum $1/p+1/q+1/r$ less than~1 but as 
close to~1 as possible. The triple $(2,3,7)$ obviously provides the answer. Then
$$
2-2g \eq n\left(\frac{1}{2}+\frac{1}{3}+\frac{1}{7}-1\right) \eq -\frac{n}{42}
$$
which gives $n=84(g-1)$.

\begin{remark}
For convenience of drawing, we will regard $\Delta$ as the triangle group 
$\Delta(3,2,7)$, so that $\R$ and its quotient dessins $\D$ have
type $(3,2,7)$ rather than $(2,3,7)$. Writing the periods in that 
order means that the black and white vertices have valencies 
dividing $3$ and $2$ respectively; this allows us in Figures~4, 5 and 16 to represent dessins more 
simply as uncoloured maps, by leaving the white vertices implicit 
(so that those of valency $1$ give rise to free edges). However, in Figures 7, 9 and 10 only white vertices of valency $2$ are omitted, and elsewhere we will show all white vertices explicitly.
\end{remark}

Many non-abelian finite simple groups are now known to be (or known not to be) 
Hurwitz groups: see the surveys by Conder~\cite{Conder-90, Conder-10}. For example, 
we have the following theorem~\cite{Macbeath-69}:

\begin{theorem}[Macbeath]\label{Macbthm}
The group $G = {\rm PSL}_2(q)$ is a Hurwitz group if and only if one of the
following holds:
\begin{itemize}
\item[\rm 1.] $q = 7$, or
\item[\rm 2.] $q$ is a prime $p\equiv \pm 1$ \mbox{\rm mod } $(7)$, or
\item[\rm 3.] $q = p^3$ for some prime $p\equiv \pm 2$ or $\pm 3$ \mbox{\rm mod } $(7)$.
\end{itemize}
In cases $(1)$ and $(3)$ the Hurwitz surface and dessin associated with $G$ are unique, 
but in case $(2)$ there are three of each, corresponding to three normal subgroups $N$ 
of $\Delta$ with quotient $G$.
\end{theorem}

Of course, by Dirichlet's theorem on primes in arithmetic progressions there are infinitely many examples satisfying each of the congruences in cases (2) and (3).

%%%%%%%%%%%%%%%%%%

\section{Klein's curve}

In this section we take $\R$ to be the regular dessin of type $(3,2,7)$ on Klein's 
quartic curve
\[u^3v+v^3w+w^3u=0,\]
the Hurwitz surface of least genus, namely $g=3$. The automorphism group $G$ of 
both $\R$ and the curve is the smallest Hurwitz group, namely 
${\rm GL}_3(2)={\rm SL}_3(2)={\rm PSL}_3(2)$, of order $168$, also isomorphic 
to ${\rm PSL}_2(7)$. We will try to represent $\R$ as the minimal regular cover 
$\widetilde\D$ of a smaller dessin $\D = \R/H$ where $H\le G$.

%%%%%%%%

\subsection{${\rm PSL}_3(2)$}

Since $G$ is simple, it acts faithfully on the cosets of any proper subgroup. The subgroups of $G$ are well-known, and the smallest index of any proper subgroup $H$ is $7$, with $H$, the stabiliser of a point in the Fano plane $\Pi:=\P^2(2)$, isomorphic to the symmetric group ${\rm S}_4$. We therefore first look for quotient dessins $\D=\R/H$ of degree $7$, arising from the action of $G$ as the automorphism group of $\Pi$.

\begin{figure}[htbp]
\begin{center}
\begin{tabular}[b]{ccc}
0 & 1 & 3 \\
1 & 2 & 4 \\
2 & 3 & 5 \\
3 & 4 & 6 \\
4 & 5 & 0 \\
5 & 6 & 1 \\
6 & 0 & 2 \\
\\
\\

\end{tabular}
\hspace{2cm}
\includegraphics[scale=0.5]{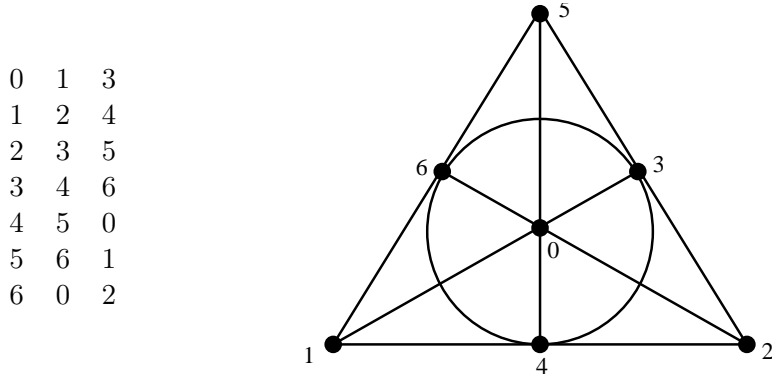}
\end{center}
\caption{The Fano plane.}
\label{fig:fano}
\end{figure}

Let us construct $\Pi$ by using the {\em difference set}\/ $\{0,1,3\}$ (equivalently, 
the set $\{1, 2, 4\}$ of quadratic residues) in the additive group of the field $\F_7$: 
every non-zero residue modulo 7 appears exactly once as the difference of two elements of this set. Then the lines are the translates of this set, see Figure~\ref{fig:fano}.

\begin{figure}[htbp]
\begin{center}
\includegraphics[scale=0.4]{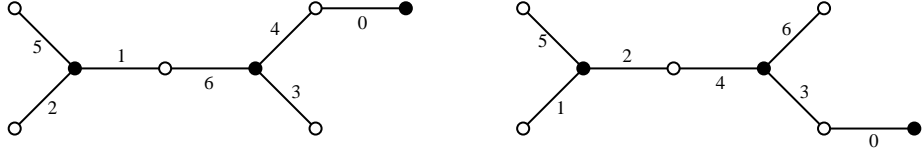}
\end{center}
\caption{\small Trees providing a $(3,2,7)$-generation of the group ${\rm PSL}_3(2)$.}
\label{fig:trees-labeled}
\end{figure}

By construction, this geometry has the following automorphism:
$$
\gamma \eq (0,1,2,3,4,5,6): a \mapsto a+1 \,\, {\rm mod}\,\, (7).
$$
Consider a clockwise rotation of our figure through $120^{\circ}$. The point 0
is fixed, so we get the permutation $x$ which we will use as a black vertex
permutation of a dessin we are looking for:
$$
x \eq (1,5,2)(3,4,6).
$$
This permutation gives us the left tree in Figure \ref{fig:trees-labeled}.
The permutation $z$ defining the (only) face in this case is equal to 
$$
z \eq \gamma^3 \eq (0,3,6,2,5,1,4): a \mapsto a+3 \,\, {\rm mod}\,\, (7).
$$
The tree on the right in the same figure is given by permutations 
$x^{-1}$ and $z^{-1}=\gamma^4$.

The involutions which correspond to the white vertices may be interpreted as 
mirror symmetries of our geometry with respect to a line. 
Indeed, the involution $y = (0,4)(1,6)$, corresponding to the left tree, 
preserves point-wise the line $\{2,3,5\}$ and preserves set-wise the lines 
$\{5,0,4\}$ and $\{5,1,6\}$. Hence, it is the ``mirror symmetry'' with respect 
to the line $\{2,3,5\}$. Analogously, the involution of the right tree is the 
symmetry with respect to the line $\{5,6,1\}$.

\begin{remark}
We did not prove that the monodromy group $G$ of the trees is indeed the whole 
group ${\rm PSL}_3(2)$ and not a proper subgroup. However, $G$ clearly 
has order divisible by $3\cdot 2\cdot7=42$, so it has index $m\le 4$ in
${\rm PSL}_3(2)$. Hence its core (intersection of all 
conjugates) is a normal subgroup of index at most $4!$, since it is the kernel of the homomorphism ${\rm PSL}_3(2)\to{\rm S}_m$ induced by the action of ${\rm PSL}_3(2)$ on the cosets of $G$. 
Since ${\rm PSL}_3(2)$ is simple and of order greater than $4!$, the core, and hence also $G$, must be the whole group.
\end{remark}

\begin{remark}
It is interesting to note that Klein's curve is a remarkable (and complicated) object. 
A book \cite{Levy-99} of 340 pages has been devoted to the study of various properties 
of this curve. However, our simple picture {\em proves the existence}\/ of such a curve.
\end{remark}

%%%%%%%%

\subsection{Character table and Frobenius's formula for ${\rm PSL}_3(2)$}

The above construction is consistent with the information provided to us by the
character table of the group ${\rm PSL}_3(2)$. The table, computed by GAP, is 
shown in Table~\ref{char-table-psl_3_2}. (It can also be found
in the {\sc Atlas}~\cite{Atlas-05}, where this group is denoted by $L_2(7)$.)
 
     The six conjugacy classes of elements of the group are denoted in the table by 
     {\tt 1a,\,2a}, etc. The notation used in the {\sc Atlas}, which we adopt 
     below, is $1A,2A,\ldots$ The first digit gives the order of the elements of
     the class. We see, for example, that there are two conjugacy classes $7A$ and $7B$ of
     elements of order~7.
     The notation {\tt X.1,\,X.2}, etc. is used for the irreducible characters; we will 
     use more standard {\sc Atlas} notation $\chi_1,\chi_2,\ldots$ for characters. 
     A~dot in the table means zero; {\tt E(7)} is the primitive 7th root of 
     unity; finally, {\tt /A} means the complex conjugate of {\tt A}.

The character of the permutation representation of degree~7 is $\chi_1+\chi_4$; 
its value on a given class is the number of fixed points. We see that for the 
class $2A$ it is equal to 3, thus giving us the cycle-length partition $2^21^3$; 
for $3A$ this value is 1, which gives us the partition $3^21^1$; and for the 
classes $7A$ and $7B$ the number of fixed points is, naturally, zero.

\begin{table}[htbp]
\begin{center}
\begin{minipage}{5.4cm}
\begin{verbatim}
       1a 2a 4a 3a 7a 7b

X.1     1  1  1  1  1  1
X.2     3 -1  1  .  A /A
X.3     3 -1  1  . /A  A
X.4     6  2  .  . -1 -1
X.5     7 -1 -1  1  .  .
X.6     8  .  . -1  1  1

A = E(7)^3+E(7)^5+E(7)^6
  = (-1-Sqrt(-7))/2 = -1-b7
\end{verbatim}
\end{minipage}

\caption{\small Character table of the group ${\rm PSL}_3(2)\cong{\rm PSL}_2(7)$
computed by GAP.}
\label{char-table-psl_3_2}
\end{center}
\end{table}

The sizes of the classes $3A$, $2A$ and $7A$, also computed by GAP, are, respectively,
56, 21 and 24. Therefore, the Frobenius formula for the passport $(3A,2A,7A)$
gives us 168 triples $(x, y, z)$ with entries in these classes and $xyz=1$. 
We obtain the same number for the passport $(3A, 2A, 7B)$, so the total number 
of triples of type $(3, 2, 7)$ is $336=2|G|$. As shown in Remark~3.1 they 
all generate $G$. Now ${\rm Aut}(G)$ is an extension of $G$ by a cyclic group ${\rm C}_2$, induced by duality of the plane $\Pi$. (Alternatively, if we identify $G$ with ${\rm PSL}_2(7)$ as in Section~\ref{PSL_2(7)}, then ${\rm Aut}(G)$ can be identified with ${\rm PGL}_2(7)$; see Section~\ref{sec:quotients} for more details.) Dividing by $|{\rm Aut}(G)|=2|G|$ we see that there is 
one regular dessin $\R$ of type $(3,2,7)$ with automorphism group $G$. 

Finally, the presence of $\sqrt{-7}$ in the character table suggests that the 
trees, being considered as {\em dessins d'enfants}, should be defined over the field 
$\Q(\sqrt{-7})$. And, indeed, the computation of Bely{\u\i} functions confirms 
this hypothesis:
$$
f(t) \eq K \cdot (t^2+7a)^3\,(t-7),
$$
$$
f(t)-1 \eq K \cdot (t^2-6t+a)^2\,(t^3+5t^2+(19a+24)t+(83a+108)),
$$
where
$$
K \eq - \frac{1}{2^6 3^3 (7a+17)},
$$
and $a$ is a root of
$$
a^2+3a+4,
$$
that is,
$$
a \eq -\frac{3}{2} \pm \frac{1}{2}\sqrt{-7}.
$$
Choosing one of the two values of $a$ we get one of the two trees. Namely, 
its two black vertices of valency~3 are roots of $t^2+7a$; its only black vertex 
of valency~1 lies at the point $t=7$; its two white vertices of valency~2 are roots 
of $t^2-6t+a$; and, finally, its three white vertices of valency~1 are roots of
$t^3+5t^2+(19a+24)t+(83a+108)$. Being a polynomial, $f$ has a pole of multiplicity~7
at infinity: it corresponds to the (only) face of valency~7. For any 
$s \notin \{0,1,\infty\}$ the equation $f(t)=s$ does not have multiple roots.

\begin{figure}[htbp]
\begin{center}
\includegraphics[scale=0.4]{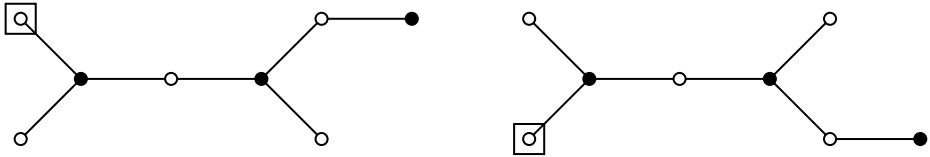}
\end{center}
\caption{\small White vertex of valency 1 which is not Galois-conjugate to the
other two.}
\label{fig:separate-vertex}
\end{figure}

An interesting observation is that the cubic factor
$$
P(t) \eq t^3+5t^2+(19a+24)t+(83a+108)
$$
in the function $f-1$ factorizes over the field $\Q(\sqrt{-7})$. Namely, 
$P = Q \cdot R$ where
$$
Q(t) = t^2 + (1+\sqrt{-7})t + \frac{-31+13\sqrt{-7}}{2}, \qquad R(t) = t + (4-\sqrt{-7}).
$$
Thus, one of the three white vertices of valency 1 is separated
 from the other two: it does not belong to the same Galois orbit.
This vertex is shown in Figure \ref{fig:separate-vertex}.

%%%%%%%%

\subsection{${\rm PSL}_2(7)$}\label{PSL_2(7)}

\begin{figure}[htbp]
\begin{center}
\includegraphics[scale=0.5]{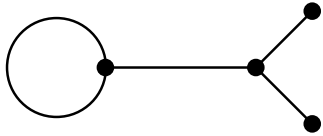}
\caption{\small A map with the monodoromy group ${\rm PSL}_2(7)$. White vertices
in the middle of the edges are implicit.}
\label{fig:psl_2_7}
\end{center}
\end{figure}

Ever since the time of Galois it has been known that ${\rm PSL}_3(2)$ is 
isomorphic (as an abstract group) to ${\rm PSL}_2(7)$: the former group acts 
on the seven points of the Fano plane $\Pi$, whereas the latter acts on the 
eight points of the projective line $\P^1(7)$ over the field $\F_7$, that is, 
on the set
$$
\{0,1,2,3,4,5,6,\infty\}.
$$
The character of this permutation representation is $\chi_1+\chi_5$, see 
Table~\ref{char-table-psl_3_2}. Therefore, the permutations in the class 
$3A$ have two fixed points, those in the class $2A$ have none, and for the 
classes $7A$ and $7B$ the number of fixed points is one. We thus obtain the 
passport $(3^2 1^2, 2^4, 7^1 1^1)$.
It is easy to see that there exists only one map with this passport: it is shown
in Figure~\ref{fig:psl_2_7}. Once again, we can only marvel at the fact that this
simple dessin provides us a $(3,2,7)$-presentation of the group ${\rm PSL}_2(7)$
and thus, in the bud, with all the information we need in order to construct the 
Klein curve.

\begin{figure}[htbp]
\begin{center}
\includegraphics[scale=0.5]{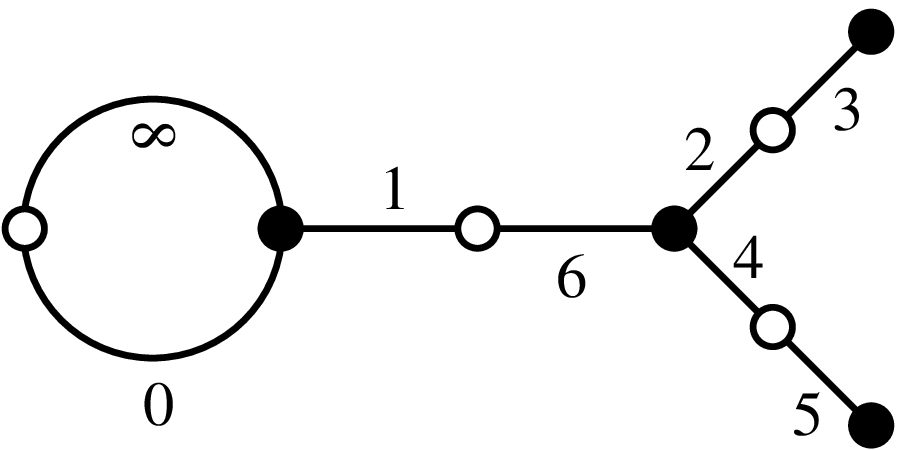}
\hspace{1cm}
\includegraphics[scale=0.5]{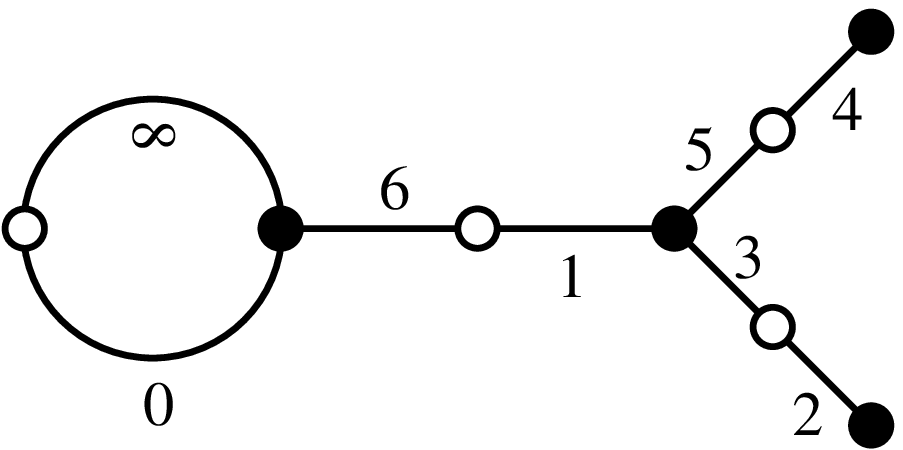}
\caption{\small $(3,2,7)$-generators of the group ${\rm PSL}_2(7)$. In the right
map, as compared with the left one, the permutation $z \in 7A$ is replaced with
$z^{-1} \in 7B$. Notice that the involution $y$ remains the same.}
\label{fig:psl_2_7-labeled}
\end{center}
\end{figure}

Labels shown in Figure \ref{fig:psl_2_7-labeled} (left) give the following
permutations obviously belonging to ${\rm PSL}_2(7)$: 
$$
y: a \mapsto -\frac{1}{a}, \qquad z: a \mapsto a+1,
$$
while $x$ can be easily computed using the relation $xyz=1$:
$$
x \eq (yz)^{-1}: a \mapsto -\frac{1}{a-1}.
$$

\begin{remark}
We multiply permutations from left to right, as is usual in symbolic calculation 
systems. But if we want, in our case, to represent the above fractional linear 
functions $x,y,z$ as $2\times 2$-matrices $X,Y,Z$, then these matrices should be 
multiplied in the inverse order: $xyz=1$ but $ZYX=I$.
\end{remark}

%%%%%%%

\subsection{Uniqueness of quotients}
\label{sec:quotients}

Why, in the case of ${\rm PSL}_2(7)$, is the quotient dessin unique (and defined 
over $\Q$), whereas in the case of ${\rm PSL}_3(2)$ we obtain a chiral pair defined 
over $\Q(\sqrt{-7})$? The answer is that $G$ has only one conjugacy class of 
subgroups $H$ of index $8$, but two of index $7$, and non-conjugate subgroups 
give non-isomorphic quotients. The subgroups of index $7$ are the stabilisers 
of points or of lines in the Fano plane $\Pi$, forming two conjugacy classes 
of size $7$, all isomorphic to ${\rm S}_4$. These two classes are transposed 
by the outer automorphism group ${\rm Out}(G)$ (corresponding to the duality of $\Pi$), just as this group 
transposes the conjugacy classes $7A$ and $7B$ of elements of order $7$ in $G$. 
We therefore have four possible choices when forming quotient dessins of degree $7$: 
there are two possibilities (up to conjugacy) for the subgroup $H$, and two (again 
up to conjugacy) for the element $z$. Since ${\rm Out}(G)$ transposes the two 
choices in each case, it has two orbits on these four possibilities, corresponding 
to two non-isomorphic dessins, transposed by changing our choice for $H$ or for $z$ 
(but not both!). In the case of quotients of degree $8$, however, there is a single 
conjugacy class of subgroups $H$ of this index, namely the stabilisers of points 
in $\P^1(7)$; our two possible choices for the class of $z$ are transposed by 
${\rm Out}(G)$, so they lead to isomorphic dessins.

In general, the number of faithful quotients of degree $n$ of a regular map $\mathcal R$ is the number of conjugacy classes of subgroups of index $n$ in $G={\rm Aut}({\mathcal R})$ with trivial core. Of course, when $n=|G|$ there is a unique quotient, namely $\mathcal R$ itself, corresponding to the identity subgroup.

\begin{figure}[htbp]
\begin{center}
%\epsfig{file=Figures/kleinmap.eps,width=5cm}
%\hspace{1cm}
%\epsfig{file=Figures/Klein.eps,width=6cm}
\includegraphics[scale=0.35]{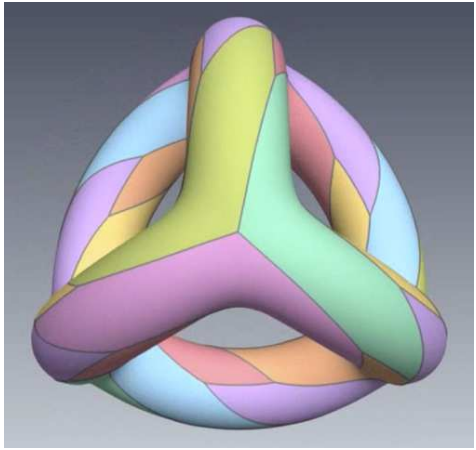}
\caption{\small The Klein map of genus 3. The author of the picture is Jarke van Wijk.}
\label{fig:Klein}
\end{center}
\end{figure}

We finish this section with a picture of the Klein map on the surface of
genus~3. The map contains 24 heptagonal faces and 56 vertices of degree~3.
The picture, see Figure~\ref{fig:Klein}, was created by Jarke van Wijk from
the Eindhoven University of Technology.

%%%%%%%%%%%%%%%%%%%%%%%%%

\section{The Fricke--Macbeath curve} 

In 1899, in the same issue of the journal ``Mathematische Annalen'', two papers 
were published. In the first one \cite{Burnside-1899}, by Burnside, it was shown
that the simple group ${\rm PSL}_2(8)$ has a $(3,2,7)$-presentation; 
according to a tradition of that era, no motivation for this result was given. In 
the second paper~\cite{Fricke-1899}, Fricke constructed a Riemann surface of genus~7
with the automorphism group ${\rm PSL}_2(8)$ of order $504=84\cdot 6$, that is,
of maximal size for that genus. Later on this surface was rediscovered by 
Macbeath \cite{Macbeath-65}.

\begin{figure}[htbp]
\begin{center}
\includegraphics[scale=0.5]{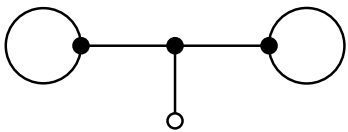}
\caption{\small A map with monodromy group ${\rm PSL}_2(8)$.}
\label{fig:psl_2_8}
\end{center}
\end{figure}

Thus, for us, in order to prove the result of Burnside and Fricke, it suffices
to produce a $(3,2,7)$-presentation of the group ${\rm PSL}_2(8)$. Such a presentation 
is shown in Figure~\ref{fig:psl_2_8}, by which we mean that this map of type $(3,2,7)$ 
has monodromy group ${\rm PSL}_2(8)$.

It would not be an easy task to give a ``purely mathematical'' proof of the above 
statement, that is, one which is not computer-assisted. Indeed, since this dessin 
of degree $9$ corresponds to the natural representation of $G$, the edge labels should 
be the elements of the projective line $\P^1(8)$ over $\F_8$, and the elements of 
this field should, in turn, be represented as polynomials of degree~2 with coefficients 
in $\F_2$. But today, using Maple or GAP, we can in a fraction of a second compute 
the order of the group, which is 504, and then look at the catalogue \cite{ButMcK-83} 
and see that there is only one permutation group of degree 9 and of order 504.

\begin{remark}
An ambitious question we would like to ask is as follows: what kind of information
about a ``big'' regular map with an automorphism group~$G$ can we extract from a 
``small'' quotient map with a monodromy group~$G$? In the above example, it is easy 
to see that Figure~\ref{fig:psl_2_8} is the only possible $(3,2,7)$-map of degree~9.
Therefore, it is defined over $\Q$. We may (or may we, indeed?) infer from this
that a $(3,2,7)$-presentation of the group ${\rm PSL}_2(8)$ is unique up to an
automorphism, and therefore the Hurwitz map of genus 7 is also unique and defined 
over $\Q$. However, Manfred Streit \cite{Streit-00} showed that the latter map
cannot be {\em realized}\/ over $\Q$, and Rub\'en Hidalgo~\cite{Hidalgo-15} 
managed to realize it over $\Q(\sqrt{-7})$. Can this be seen from 
Figure~\ref{fig:psl_2_8}?
\end{remark}

\begin{figure}[htbp]
\begin{center}
\includegraphics[scale=0.4]{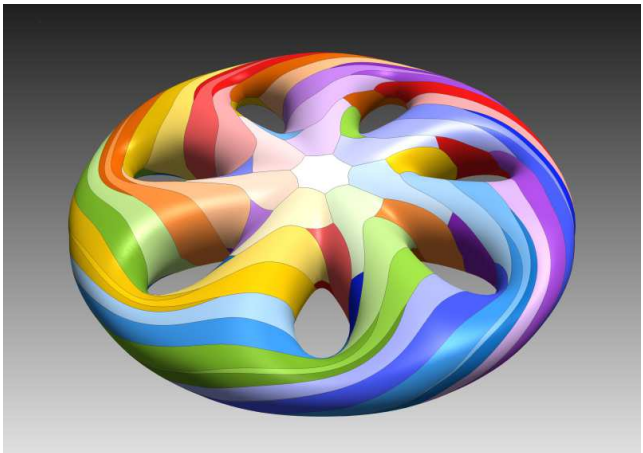}
\caption{\small The Fricke--Macbeath map of genus 7 with the automorphism group 
${\rm PSL}_2(8)$ of order $504=84\cdot 6$. The author of the picture is Jarke van Wijk.}
\label{fig:Fricke-Macbeath}
\end{center}
\end{figure}

Figure \ref{fig:Fricke-Macbeath} shows the Fricke--Macbeath map. This chef-d'oeuvre
of computer graphics was made by Jarke van Wijk~\cite{vanWijk-14}. According to van Wijk, the pictures of this kind are produced by {\em tubification}\/ of graphs~\cite{vanWijk-16}. We take a graph, chosen to have Betti number $E-V+1$ equal to the genus of the map, where $V$ and $E$ are the numbers of vertices and edges of the graph. This graph is then embedded in three-dimensional space, if possible exhibiting some group of symmetries which it shares with the map. 
We then replace its vertices and edges with spheres and tubes, thus creating a surface 
on which the desired map may be drawn. The main difficulty of this approach is that 
there are infinitely many non-isomorphic graphs with a given Betti number (and not all 
of them planar, by the way). Permutations representing the map give us complete 
information about the map itself, but tell us absolutely nothing about a convenient 
structure of a graph to be tubified and about its possible embedding into ${\mathbb R}^3$.

Figures \ref{fig:trees-labeled}, \ref{fig:psl_2_7} and \ref{fig:psl_2_8}
constitute a complete list of $(3,2,7)$-maps with a single face of degree~7.
In subsequent sections we will examine all the $(3,2,7)$-maps with 
two faces of degree~7.

%%%%%%%%%%%%%%%%%%%%%%%%%%

\section{Three Hurwitz maps of genus 14}
\label{sec:genus-14}

For our next example we take $G$ to be the Hurwitz group ${\rm PSL}_2(13)$ of genus~$14$ 
and order $84\cdot 13 = 1092$. As for all groups ${\rm PSL}_2(q)$ with $q>11$, 
the transitive permutation representation of least degree is the natural representation, 
of degree $q+1 = 14$, so we look for possible quotient dessins $\D$ of this degree. 
The point stabilisers are isomorphic to ${\rm C}_{13}\rtimes{\rm C}_6$. The group $G$ has unique conjugacy classes $3A$ and $2A$ of elements of order $3$ and $2$, 
but it has three self-inverse conjugacy classes, $7A$, $7B$ and $7C$, of elements of 
order $7$; those in $7B$ and $7C$ are the squares and fourth powers of those in $7A$. 
(This applies to all the Hurwitz groups ${\rm PSL}_2(q)$ with $q\ne 7$.) For each 
of these classes, if we combine it with the classes $3A$ and $2A$, the Frobenius 
formula gives $2|G|$ triples, all generating $G$ since no proper subgroup of $G$ 
is a Hurwitz group. (There are four conjugacy classes of maximal subgroups of $G$, isomorphic to ${\rm C}_{13}\rtimes{\rm C}_6$, ${\rm D}_7$, ${\rm D}_6$ and ${\rm A}_4$; none of them has order divisible by 2, 3 and 7.) Since ${\rm Aut}(G)={\rm PGL}_2(13)$ has order $2|G|$, 
leaving each of the three classes of elements of order $7$ invariant,
we obtain three regular dessins $\mathcal R$, one for each class. As 
Streit~\cite{Streit-00} has shown, they form a Galois orbit, defined 
over the field $\Q(\cos(2\pi/7))$. (The entries in the character table of 
$G$ also belong to this field.)

\begin{figure}[htbp]
\begin{center}
\includegraphics[scale=0.5]{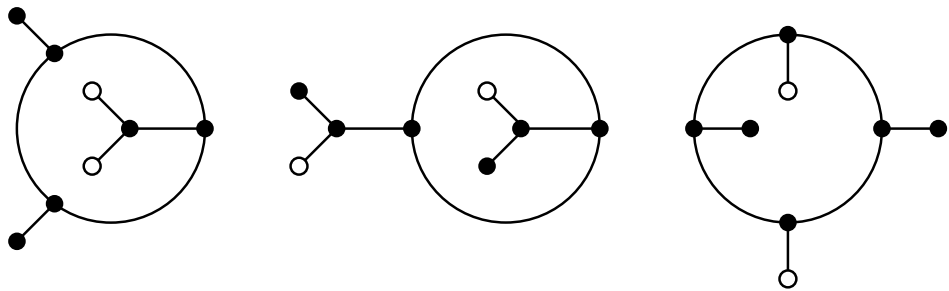}
\caption{\small Hurwitz generation of ${\rm PSL}_2(13)$.}
\label{fig:psl_2_13}
\end{center}
\end{figure}

The corresponding quotient dessins $\D$ of degree $14$ are shown in 
Figure~\ref{fig:psl_2_13}. They also form a single Galois orbit, so at first it 
is a little disturbing to see that the map on the left is invariant under a 
reflection, whereas the other two are not, and nor are they mirror images of 
each other. However, if we draw these two
as maps on the sphere, then in each case the antipodal involution provides the 
expected orientation-reversing isomorphism. These two maps also illustrate the 
Couveignes--Filimonenkov phenomenon: they are defined, but cannot be realized 
over a real field. Does their antipodal symmetry
mean something interesting for the corresponding regular maps $\R$ of genus~14?

This example illustrates a general phenomenon, in which all the Mac\-beath--Hurwitz 
groups ${\rm PSL}_2(q)$ for $q\ne 7$ or $27$ have unique classes of elements of orders 
$3$ and $2$, but three of order $7$; one can show by direct calculation with M\"obius transformations or by the Frobenius formula that each choice of classes of elements of these orders gives $2|G|$ generating triples of type $(3,2,7)$, resulting in a total of $6|G|$ triples. 

Here we obtained three Hurwitz dessins $\mathcal R$, one for each class of elements of order 7, so why did this not happen in our earlier example, where $G = {\rm PSL}_2(8)$? The answer is that in that case, and indeed in all examples of case (3) of Theorem~\ref{Macbthm}, where $q=p^3$ for some prime $p$, ${\rm Aut}(G)$ is not ${\rm PGL}_2(q)$ but the larger group ${\rm P\Gamma L}_2(q)$. This is an extension of ${\rm PGL}_2(q)$ by a cyclic group of order $3$ induced by the Galois group ${\rm Gal}(\F_q)\cong {\rm C}_3$ of the field $\F_q$. It acts by permuting the three classes of elements of order $7$ transitively, so that the $6|G|$ triples form a single orbit under ${\rm Aut}(G)$ and hence correspond to a single Hurwitz dessin. In case (2) of Theorem~\ref{Macbthm}, with $q=p$ prime, the Galois group of the field is trivial, so the triples form three orbits, corresponding to three Hurwitz dessins.

Pictures of the corresponding regular maps of genus 14 are not yet available, certainly because of the problem with their spatial arrangement.

%%%%%%%%%%%%%%%%%%%%%%%

\section{Genus 17}\label{genus17}

\begin{figure}[htbp]
\begin{center}
\includegraphics[scale=0.3]{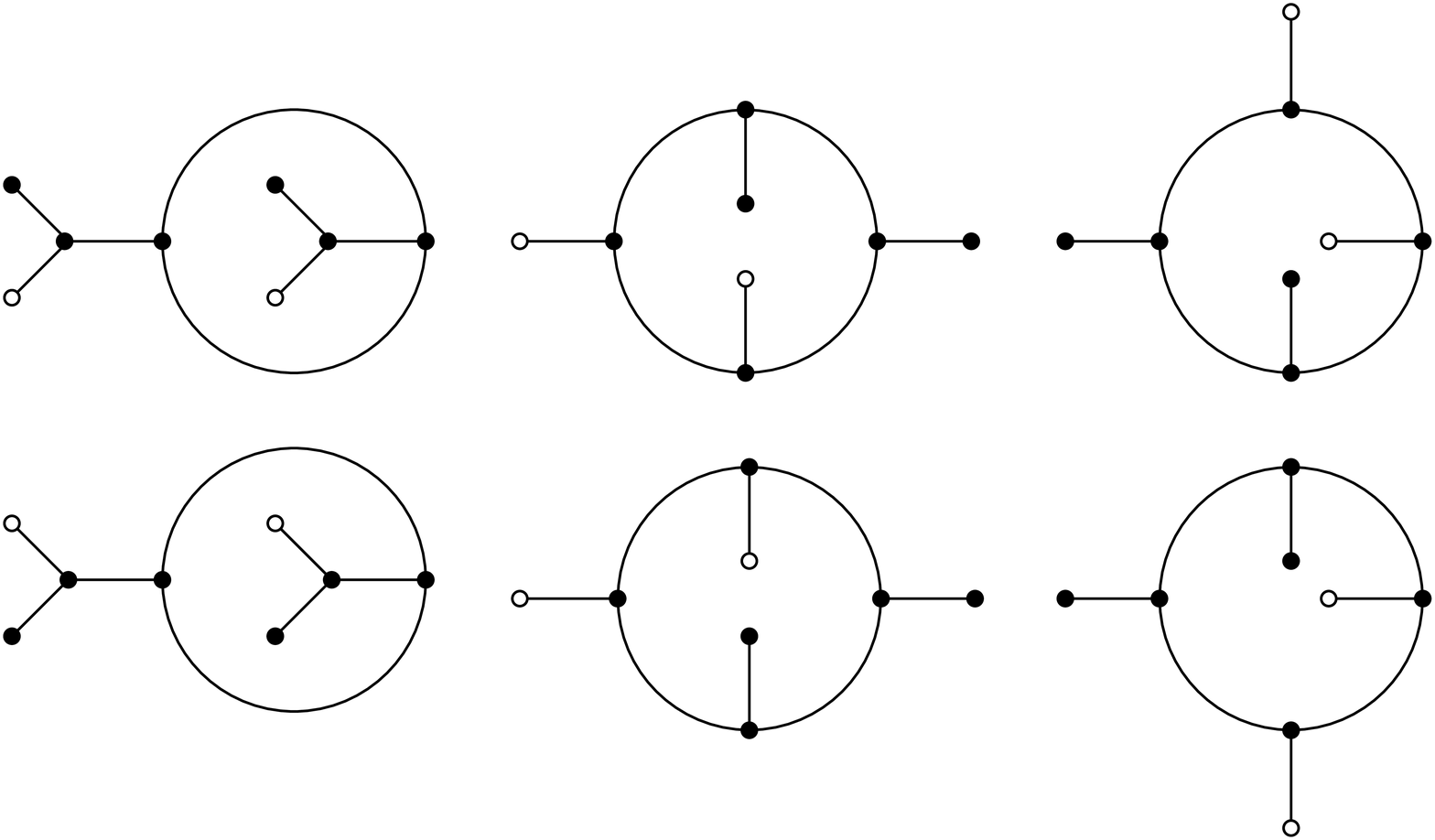}
\caption{\small Six imprimitive maps. There are seven blocks of size 2.}
\label{fig:asl_3_2}
\end{center}
\end{figure}

There remain nine $(3,2,7)$-maps with two faces of degree 7. Six of them are
shown in Figure~\ref{fig:asl_3_2}, three more are given in Figure~\ref{fig:a_15}.

\begin{figure}[htbp]
\begin{center}
\includegraphics[scale=0.3]{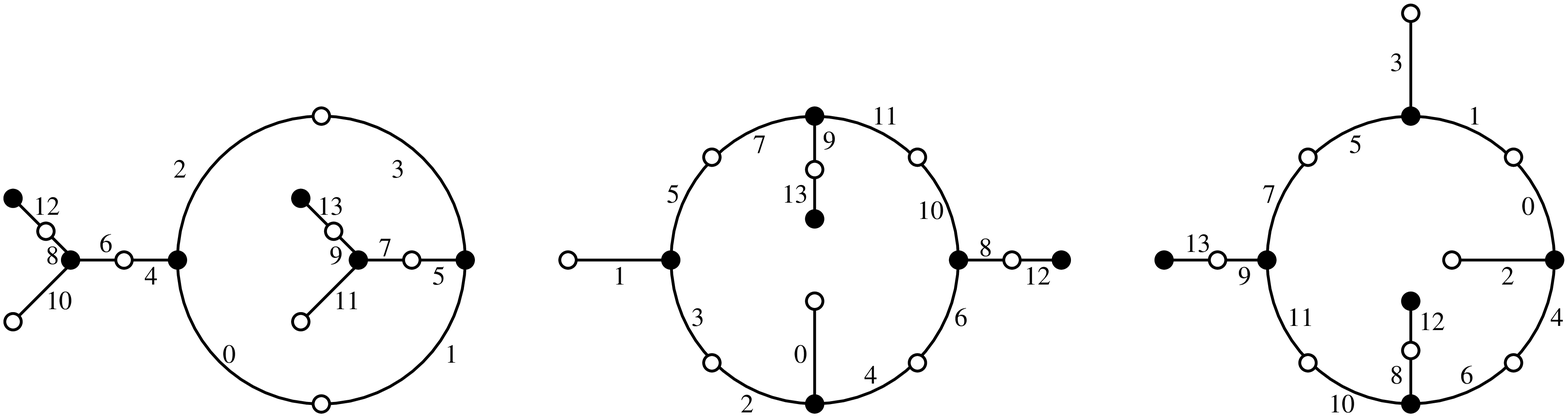}
\caption{\small Labelling of the three upper maps in Figure \ref{fig:asl_3_2}.}
\label{fig:asl_3_2-labeled}
\end{center}
\end{figure}

\begin{figure}[htbp]
\begin{center}
\includegraphics[scale=0.25]{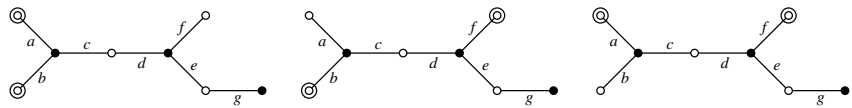}
\caption{\small ${\rm PSL}_3(2)$-trees, labelled with blocks for the maps in Figure 
\ref{fig:asl_3_2-labeled}.}
\label{fig:asl_3_2-blocks}
\end{center}
\end{figure}

All six maps in Figure~\ref{fig:asl_3_2} are imprimitive: they cover the 
${\rm PSL}_3(2)$-trees we have seen before. Figure~\ref{fig:asl_3_2-labeled} 
proposes a labelling of the maps in the first row of Figure~\ref{fig:asl_3_2} for which 
the blocks are the same. These blocks are shown in Table~\ref{tab:asl_3_2-blocks},
while the action of the permutations on blocks and the ramification points may
be seen in Figure~\ref{fig:asl_3_2-blocks}.

\begin{table}[htbp]
\begin{center}
\begin{tabular}{|c|c|c|c|c|c|c|}
\hline
$a$ & $b$ & $c$ & $d$ & $e$ & $f$  & $g$ \\
\hline\hline
0   & 2   & 4   & 6   & 8   & 10   & 12 \\
\hline
1   & 3   & 5   & 7   & 9   & 11   & 13 \\
\hline
\end{tabular}
\caption{Blocks of the permutations describing maps in Figure \ref{fig:asl_3_2-labeled}.}
\label{tab:asl_3_2-blocks}
\end{center}
\end{table}

Indeed, let us write the permutations describing the maps of 
Figure~\ref{fig:asl_3_2-labeled} and their action on the blocks. 
Here are permutations which correspond to the left map, and also to the 
left tree of Figure~\ref{fig:asl_3_2-blocks}.
We see that the permutations $x$ and $z$ are unramified while $y$ is ramified
over the vertices $a$ and $b$ of the tree. The corresponding cycles of $y$ are
underlined and the vertices of the tree are singled out. It is clear that the map
which is mirror symmetric to this one will be ramified over the mirror symmetric tree.
\begin{eqnarray*}
x & = & (0,2,4)(1,3,5)(6,8,10)(7,9,11)(12)(13) \\
  & = & (a,b,c)(a,b,c)(d,e,f)(d,e,f)(g)(g), \\
y & = & (0,1)(2,3)(4,6)(5,7)(8,12)(9,13)(10)(11) \\
  & = & (\underline{a,a})(\underline{b,b})(c,d)(c,d)(e,g)(e,g)(f)(f), \\
z & = & (0,5,11,9,13,7,3)(1,4,10,8,12,6,2) \\
  & = & (a,c,f,e,g,d,b)(a,c,f,e,g,d,b).
\end{eqnarray*}
Permutations for the map in the middle: they are ramified over $b$ and $f$ and thus
correspond to the tree in the middle of Figure~\ref{fig:asl_3_2-blocks}.
\begin{eqnarray*}
x & = & (0,2,4)(1,3,5)(6,8,10)(7,9,11)(12)(13) \\
  & = & (a,b,c)(a,b,c)(d,e,f)(d,e,f),(g)(g) \\
y & = & (0)(1)(2,3)(4,6)(5,7)(8,12)(9,13)(10,11) \\
  & = & (a)(a)(\underline{b,b})(c,d)(c,d)(e,g)(e,g)(\underline{f,f}), \\
z & = & (0,4,10,9,13,7,3)(1,5,11,8,12,6,2) \\
  & = & (a,c,f,e,g,d,b)(a,c,f,e,g,d,b).
\end{eqnarray*}
Finally, the permutations for the map on the right: they are ramified over $a$ and $f$.
\begin{eqnarray*}
x & = & (0,2,4)(1,3,5)(6,8,10)(7,9,11)(12)(13) \\
  & = & (a,b,c)(a,b,c)(d,e,f)(d,e,f),(g)(g) \\
y & = & (0,1)(2)(3)(4,6)(5,7)(8,12)(9,13)(10,11) \\
  & = & (\underline{a,a})(b)(b)(c,d)(c,d)(e,g)(e,g)(\underline{f,f}), \\
z & = & (0,4,10,9,13,7,3)(1,5,11,8,12,6,2) \\
  & = & (a,c,f,e,g,d,b)(a,c,f,e,g,d,b).
\end{eqnarray*}

Recall that the three white vertices of Figure~\ref{fig:asl_3_2-blocks} split 
into two Galois orbits: the vertex $b$ is separated from the other two (see
Figure~\ref{fig:separate-vertex}). Therefore we may suppose that the third map,
for which the ramification points avoid the vertex $b$, might not behave in
the same way as the other two maps. And, indeed, computing the order of the
monodromy group for all three maps we find out that it is equal to 1344
for the first and second map while it is 168 for the third. Of course, the
same is true for their mirror images.

The group of size 168 can only be ${\rm PSL}_3(2)$. The identification of the
Hurwitz group $G$ of size 1344 is a more subtle matter. The catalogue \cite{Butler-93}
shows that there are three non-isomorphic permutation groups of degree~14
and of order~1344.

One of these three groups cannot be projected onto ${\rm PSL}_3(2)$,
so it cannot arise here. This group is a semidirect product $A\rtimes B$ of 
an elementary abelian normal subgroup $A\cong ({\rm C}_2)^6$ 
by a complement $B\cong {\rm AHL}_1(7)\cong {\rm C}_7\rtimes {\rm C}_3$, the 
subgroup of index~$2$ in ${\rm AGL}_1(7)$ consisting of the affine transformations 
$t\mapsto at+b$ of $\F_7$ for which $a$ is a non-zero square. (Here `H' stands 
for `half'.) The group acts on the 14 points of the cartesian product 
$\F_2\times\F_7$, with elements of $A$ acting on pairs $(s,t)$ by preserving~$t$ 
and changing an even number of coordinates $s$, while elements of~$B$ preserve~$s$ 
while acting naturally on~$t$. This action is imprimitive, since the group permutes 
the seven pairs $\{(0,t), (1,t)\}$. Since the group is solvable, it cannot be 
a Hurwitz group. 

The other two groups are also imprimitive: each is an extension of an elementary 
abelian normal subgroup $T\cong ({\rm C}_2)^3$ by ${\rm GL}_3(2)={\rm PSL}_3(2)$, 
where $T$ is the kernel of the action on seven blocks of size $2$. The obvious 
example of such a group is the affine group ${\rm AGL}_3(2)$, the group of all 
affine transformations of a $3$-dimensional vector space $V$ over $\F_2$: 
this acts on the $14$ affine planes in $V$, with $T$ as the group of translations, 
complemented by the subgroup ${\rm GL}_3(2)$ fixing the vector $0$. However, the 
following argument shows that our Hurwitz group $G$ must be isomorphic to the 
third group, which is a {\em non-split\/} extension of $T$ by ${\rm GL}_3(2)$
while ${\rm AGL}_3(2)$ is a split extension.

If $T$ has a complement $C$ in $G$, then $C$ lifts to a subgroup $M$ of index~8 
in $\Delta=\Delta(3,2,7)$. Since $C$, being isomorphic to ${\rm GL}_3(2)$, is simple, its core in $G$ is trivial, 
so the core of~$M$ in~$\Delta$ is a normal subgroup~$N$ of index 1344 in~$\Delta$. 
This subgroup~$M$ corresponds to a dessin of degree~8 and type $(3,2,7)$. 
It is easy to see that any dessin of this degree and type must have passport 
$(3^21^2, 2^4, 7^11^1)$, and we saw in \S\ref{PSL_2(7)}
that the only possibility is the dessin in Figure \ref{fig:psl_2_7}. However, 
the monodromy group of this dessin is isomorphic to
${\rm GL}_3(2)$, so $N$ has index 168 
in $\Delta$, a contradiction. Thus $G$ is a non-split extension of $T$ 
by ${\rm GL}_3(2)$, and in particular it cannot be isomorphic to ${\rm AGL}_3(2)$, 
as is sometimes assumed. (As confirmation of this, $G$ has elements of order~$8$, such as $yz^3$, whereas ${\rm AGL}_3(2)$ does not; see also Section~\ref{17again}.)  The group thus obtained is number 33 in Butler's catalogue~\cite{Butler-93}, and 14T33 in the database~\cite{LMFDB}.

\begin{figure}[h!]

\begin{center}
\begin{tikzpicture}[scale=0.5, inner sep=0.8mm]

\node (a) at (0,0)  [shape=circle, draw, fill=black] {};
\node (b) at (-2,2) [shape=circle, draw, fill=black] {};
\node (c) at (0,4)  [shape=circle, draw, fill=black] {};
\node (d) at (2,2) [shape=circle, draw, fill=black] {};
\node (e) at (0,7) [shape=circle, draw, fill=black] {};
\node (f) at (0,-3) [shape=circle, draw, fill=black] {};

\draw (a) to (b) to (c) to (d) to (a);
\draw (c) to (e);
\draw (a) to (f);

\node at (1,7) {$\Delta$};
\node at (1,4) {$K$};
\node at (-3,2) {$N_1$};
\node at (3,2) {$N_2$};
\node at (1.5,0) {$K'K^2$};
\node at (0.8,-3) {$1$};

 \end{tikzpicture}

\end{center}
\caption{The subgroups $N_i$ of $\Delta$}
\label{N1N2}
\end{figure}
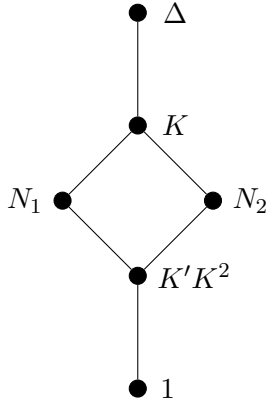

One can construct this Hurwitz group $G$ homologically as follows. The normal subgroup $K$ 
of index 168 in $\Delta$ corresponding to Klein's quartic curve is the fundamental 
group of this surface $S$ of genus 3, so its commutator quotient group $K/K'$ is the abelianised fundamental group, that is, the first integer homology group $H_1(S,\mathbb Z) \cong \mathbb Z^6$. Similarly, if $K^2$ denotes the group generated by the squares in $K$ then $K/K'K^2$ is the reduction of $H_1(S,\mathbb Z)$ mod 2, namely $H_1(S,\F_2) \cong (\F_2)^6$. This is a 6-dimensional module for the automorphism group $A \cong {\rm GL}_3(2)$ of this dessin, and it decomposes as a direct sum of two irreducible 3-dimensional $A$-submodules, corresponding to two normal subgroups $N_1$ and $N_2$ of $\Delta$ lying between $K$ and $K'K^2$ (see Figure~\ref{N1N2}). 
These subgroups $N_i$ are conjugate in the extended triangle group $\Delta[3,2,7]$, so they correspond to a chiral pair of regular dessins ${\mathcal R}_i$ of type $(3,2,7)$ and genus 17 (the duals of the chiral pair of maps C17.1 in Conder's catalogue~\cite{Conder-homepage}).
The quotient groups $\Delta/N_i$ give two realisations of the Hurwitz group $G$ as ${\rm Aut}({\mathcal R}_i)$, each having a normal subgroup $T = K/N_i \cong {\rm C}_2^3$ with quotient $\Delta/K\cong {\rm GL}_3(2)$. The minimal common cover of $\mathcal R_1$ and $\mathcal R_2$ is a Hurwitz dessin of genus 129 with automorphism group $\Delta/K'K^2$, an extension of ${\rm C}_2^6$ by ${\rm GL}_3(2)$. 

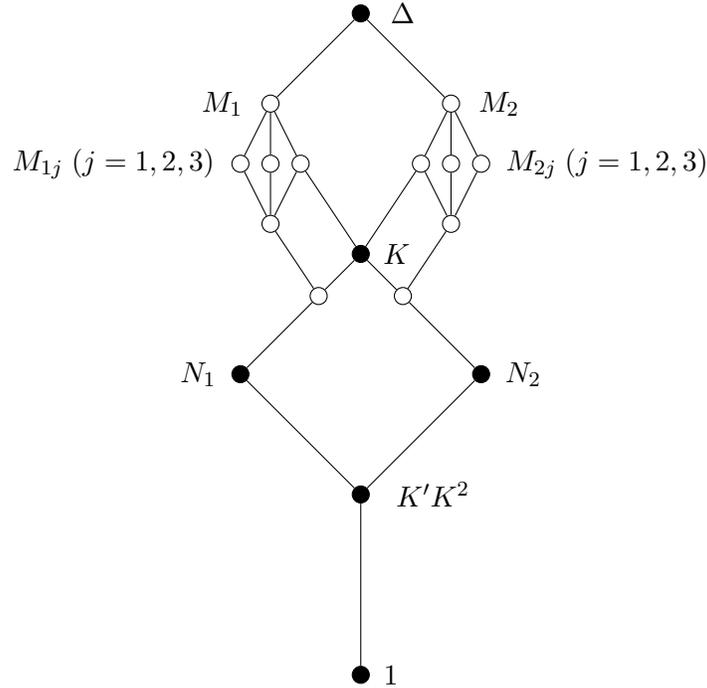
\begin{figure}[h!]

\begin{center}
\begin{tikzpicture}[scale=0.8, inner sep=0.8mm]

\node (kk) at (0,0)  [shape=circle, draw, fill=black] {};
\node (n1) at (-2,2) [shape=circle, draw, fill=black] {};
\node (k) at (0,4)  [shape=circle, draw, fill=black] {};
\node (n2) at (2,2) [shape=circle, draw, fill=black] {};
\node (d) at (0,8) [shape=circle, draw, fill=black] {};
\node (1) at (0,-3) [shape=circle, draw, fill=black] {};

\node at (0.7,8) {$\Delta$};
\node at (0.6,4) {$K$};
\node at (-2.7,2) {$N_1$};
\node at (2.7,2) {$N_2$};
\node at (1.2,0) {$K'K^2$};
\node at (-2.3,6.5) {$M_1$};
\node at (-4.1,5.5) {$M_{1j}\;(j=1,2,3)$};
\node at (2.3,6.5) {$M_2$};
\node at (4.1,5.5) {$M_{2j}\;(j=1,2,3)$};
\node at (0.5,-3) {$1$};

\node (m1) at (-1.5,6.5)  [shape=circle, draw] {};
\node (m11) at (-1,5.5)  [shape=circle, draw] {};
\node (m12) at (-1.5,5.5)  [shape=circle, draw] {};
\node (m13) at (-2,5.5)  [shape=circle, draw] {};
\node (p1) at (-1.5,4.5)  [shape=circle, draw] {};
\node (q1) at (-0.7,3.3)  [shape=circle, draw] {};

\node (m2) at (1.5,6.5)  [shape=circle, draw] {};
\node (m21) at (1,5.5)  [shape=circle, draw] {};
\node (m22) at (1.5,5.5)  [shape=circle, draw] {};
\node (m23) at (2,5.5)  [shape=circle, draw] {};
\node (p2) at (1.5,4.5)  [shape=circle, draw] {};
\node (q2) at (0.7,3.3)  [shape=circle, draw] {};

\draw (kk) to (n1) to (q1) to (k) to (q2) to (n2) to (kk);
\draw (k) to (m11) to (m1) to (d);
\draw (k) to (m21) to (m2) to (d);
\draw (kk) to (1);

\draw (m11) to (p1);
\draw (m21) to (p2);
\draw (m1) to (m12) to (p1);
\draw (m2) to (m22) to (p2);
\draw (m1) to (m13) to (p1) to (q1);
\draw (m2) to (m23) to (p2) to (q2);

 \end{tikzpicture}

\end{center}
\caption{The subgroups $M_{ij}$ of $\Delta$ (normal and non-normal subgroups 
are indicated in black and white)}
\label{Mij}
\end{figure}

The six dessins in Figure~\ref{fig:asl_3_2} correspond to six conjugacy classes of subgroups of index $14$ in $\Delta$, as follows. The two trees in Figure~\ref{fig:trees-labeled} correspond to two conjugacy classes of subgroups $H_i\cong {\rm S}_4$ in ${\rm GL}_3(2)$ for $i=1, 2$, the stabilisers of points and lines in the Fano plane. These lift to two conjugacy classes of subgroups $M_i$ of index $8$ in $\Delta$. These are Fuchsian groups whose signatures $(0; 2, 2, 2, 3)$ can be deduced from Singerman's Theorem~\cite{Singerman-70}: they have genus $0$, since this is the genus of the trees, and their elliptic periods correspond to the three fixed points of the generator $y$ of order $2$ (the three white vertices of valency $1$), and the unique fixed point of the generator $x$ of order $3$ (the black vertex of valency $1$). These groups $M_i$ therefore have presentations
\[\langle X_1, X_2, X_3, X_4\mid X_1^2=X_2^2=X_3^2=X_4^3=X_1X_2X_3X_4=1\rangle,\]
from which it is clear that they each have three subgroups $M_{ij}$ of index $2$, the normal closures in $M_i$ of $\{X_j, X_4\}$ for $j=1, 2, 3$ (see Figure~\ref{Mij}).

For each $i=1, 2$, one of the three subgroups $M_{ij}$ contains $K$, which is therefore its core: this is the lift to $\Delta$ of the commutator subgroup $H_i'\cong {\rm A}_4$ of $H_i$, giving rise to a dessin (the third in each row) with monodromy group ${\rm GL}_3(2)\cong\Delta/K$. The other two subgroups $M_{ij}$ give the first and second dessins in each row, with monodromy groups $G\cong\Delta/N_i$ of order $1344$.

%%%%%%%%%%%%

%%%%%%%%%%%%%%%%%%%%%%%%%

\section{Genus 118}
\label{sec:genus-118}

Conder \cite{Conder-87} has listed all the Hurwitz groups of order up to
one million. The next genus after 17 for which there exist Hurwitz maps is
$g=118$, with the automorphism group $G={\rm PSL}_2(27)$. As shown by Macbeath 
(see Theorem~\ref{Macbthm}), there is just one Hurwitz curve associated with 
this group.

\begin{figure}[htbp]
\begin{center}
\includegraphics[scale=0.5]{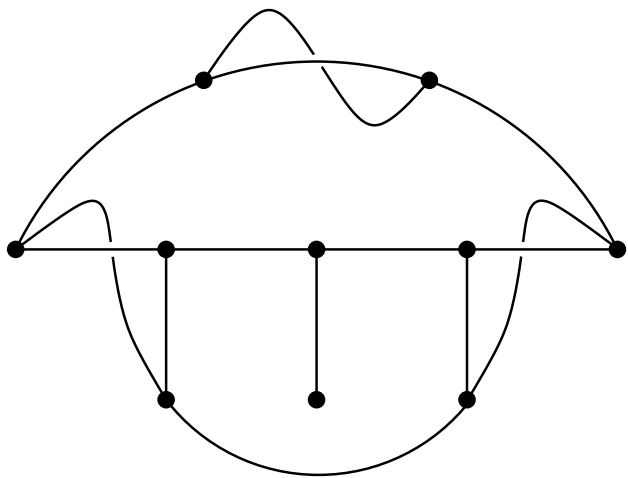}
\end{center}
\caption{\small A map of type $(3,2,7)$ with the monodromy group ${\rm PSL}_2(27)$.
This map is of genus $g=1$.}
\label{fig:psl_2_27}
\end{figure}

The natural representation of $G$ has degree $28$. The corresponding quotient map 
is shown in Figure~\ref{fig:psl_2_27}. It is drawn in the plane with three crossings, 
but it is easy to see that these can be removed by adding a single handle, so its 
genus is $1$. Alternatively, one can check that there are four faces, all of 
degree $7$, so the Euler characterisitc is $0$.

The following permutations describe the map in Figure~\ref{fig:psl_2_27} and 
generate the group ${\rm PSL}_2(27)$:
\begin{eqnarray*}
x & = & (1,2,4)(5,8,24)(6,21,10)(7,16,15)(9,25,28)(11,13,14)(12,27,23) \\ 
  &   & (17,26,18)(19,20,22), \\
y & = & (1,13)(2,25)(3,27)(4,23)(5,16)(6,12)(7,26)(8,22)(9,11)(10,17) \\ 
  &   & (14,18)(15,21)(19,24)(20,28), \\
z & = & (1,11,28,19,8,20,25)(2,9,14,26,15,6,23)(3,12,10,18,13,4,27) \\
  &   & (5,7,17,21,16,24,22).
\end{eqnarray*}
The label of the half-edge attached to the vertex of valency 1 is 3; the positions
of the other labels can easily be derived. Recall that the labels rotate around 
vertices in the counterclockwise direction.

%%%%%%%%%%%%%%%%%%%%%%%%%

\section{Higher genus examples}

\subsection{Calculating the genus}

More generally, when the group $G={\rm PSL}_2(q)$ is a Hurwitz group (as in Theorem~\ref{Macbthm}), the genus of the corresponding Hurwitz dessin (or dessins) $\R$ is
\[g=\frac{|G|}{84}+1=\frac{q(q^2-1)}{168}+1\]
for $q\ne 8$, and $g=7$ for $q=8$. Asymptotically, we thus have
\[g\sim\frac{q^3}{168}\quad{\rm as}\quad q\to\infty.\]
However, the natural representation of $G$ provides a quotient dessin $\D$ of much lower genus, which we will now calculate.

For any prime power $q>11$ the least index of any proper subgroup of ${\rm PSL}_2(q)$ is $q+1$, attained only by the conjugacy class of point-stabilisers $H$ in the natural representation on $\P^1(q)$. A non-identity M\"obius transformation (over any field) has at most two fixed points, so it follows that in this representation of $G$, any elements of orders $3$ and $7$ must have cycle-structures $3^a1^{q+1-3a}$ and $7^c1^{q+1-7c}$ where $a=\lfloor(q+1)/3\rfloor$ and $c=\lfloor(q+1)/7\rfloor$. In addition, the simplicity of $G$ implies that any element of order $2$ must induce an even permutation, so it has cycle-structure $2^b1^{q+1-2b}$ where $b=2\lfloor(q+1)/4\rfloor$.

If $G$ is a Hurwitz group, and we set aside the cases $q=7, 8$ and $27$ when $q$ is a power of $7, 2$ or $3$, we find that the numbers of black vertices, non-free edges and faces of the associated quotient map $\D=\R/H$ are therefore
\[V=\frac{q+3+2\alpha}{3},\quad E=\frac{q-\beta}{2} \quad{\rm and}\quad F=\frac{q+7+6\gamma}{7},\]
where 
\[q\equiv \alpha \, {\rm mod}\, 3,\quad q\equiv \beta \, {\rm mod}\,4
\quad{\rm and}\quad q\equiv \gamma \, {\rm mod}\, 7\]
with $\alpha, \beta, \gamma = \pm 1$. Thus the Euler characteristic of $\D$ is
\[V-E+F=2-\frac{1}{42}(q-28\alpha-21\beta-36\gamma)\]
and its genus is
\[\overline{g}=\frac{1}{84}(q-28\alpha-21\beta-36\gamma)\sim\frac{q}{84}.\]
(We have already seen that $\overline{g}=0$ when $q=7$ or $8$, and that $\overline{g}=1$ when $q=27$.)

In naive terms, whereas the information carried by the Hurwitz dessin $\R$ increases cubically with $q$, that in $\D$ increases only linearly, even though it is sufficient to determine $\R$ uniquely.

It also follows from Theorem~\ref{Macbthm} and the above calculation that the 
only cases giving planar maps $\D$ are $q=7, 8, 13, 29$ and $43$. 
Of the remaining values $q\le 100$, we obtain torus maps for $q=27, 41, 71$ and $97$, 
while for $q=83$ the genus is $2$.

(A similar calculation with the Riemann--Hurwitz formula, now not restricted to the groups ${\rm PSL}_2(q)$, shows that {\em any\/} dessin of type $(3,2,7)$ and degree $n$ has genus
\[\frac{n-28u-21v-36w}{84}+1,\]
where $u, v$ and $w$ are the numbers of fixed points of $x, y$ and $z$, that is, the numbers of black vertices, white vertices and faces which have degree $1$.)

%%%%%%%%%%%

\subsection{Failure of monotonicity}

Intuition might lead one to suppose that, among all faithful quotients of a given regular dessin, the genus should be a non-decreasing function of the degree. Indeed this follows from the Riemann--Hurwitz formula if one compares two quotients, one of which covers the other. However, there are counterexamples in which the two quotients are not comparable in this way.

For instance, let $\R$ be a Hurwitz dessin with automorphism group $G={\rm PSL}_2(q)$ 
for some prime power $q=p^e$ coprime to 2, 3 and 7
(equivalently $q\ne 7, 8$ or $27$, see Theorem~\ref{Macbthm}). If we factor out 
a Sylow $p$-subgroup $H$ of $G$, then since $|H|=q$ the resulting faithful 
quotient $\D=\R/H$ has degree $n = |G:H| = (q^2-1)/2$. By the choice of $q$, 
none of the canonical generators $x, y$ or $z$ of $G$ has fixed points in the 
action on the cosets of $H$, so $\D$ has genus
\[g = \frac{n}{84} + 1 = \frac{q^2+167}{168}.\]

However, if instead we factor out a dihedral subgroup $H'\le G$ of order $q-1$ (there is a single conjugacy class of these maximal subgroups in $G$), then the quotient dessin $\D'$ has degree $n' = q(q+1)/2 > n$. The generators $x$ and $z$ of orders 3 and 7 each lie in at most one conjugate of $H'$ (namely the normaliser of their centraliser in $G$ if $q\equiv 1$ mod 3 or mod 7), so they each have at most one fixed point in the action of $G$ on the cosets of $H'$. A simple double counting argument (also applicable to $x$ and $z$) deals with $y$: since $H' = N_G(H')$ there are $|G:H'| = q(q+1)/2$ conjugates of $H'$ in $G$, each containing $(q\pm 1)/2$ involutions as $q\equiv \pm 1$ mod 4; there is a single conjugacy class of $|G|/(q\mp 1) = q(q\pm 1)/2$ involutions in $G$, so the number of conjugates of $H'$ containing any one of them is
\[\frac{q(q+1)/2 \cdot (q\pm 1)/2}{q(q\pm 1)/2} = \frac{q+1}{2},\]
and hence this is the number of fixed points of $y$. It then follows that $\D'$ has genus
\[g' \le \frac{q(q+1)/2-21(q+1)/2}{84} + 1 = \frac{q^2-20q+147}{168} < g. \]
For example, if $q=13$, so that $n=84$ and $n'=91$, then $g=2$ whereas $g'=0$.

(The inequality $g' < g$ fails in the three excluded cases: if $q=7$ then $g=g'=0$; if $q=8$, taking $H'\cong C_7$, we have $g=0$ and $g'=1$; if $q=27$ then $g=1$ and $g'=2$.)

%%%%%%%%%%%%%%%%%%%%%%%%%%%%%%%%

\section{${\rm A}_{15}$ -- the first alternating group to arise as a \\ Hurwitz group}

\begin{figure}[htbp]
\begin{center}
\includegraphics[scale=0.4]{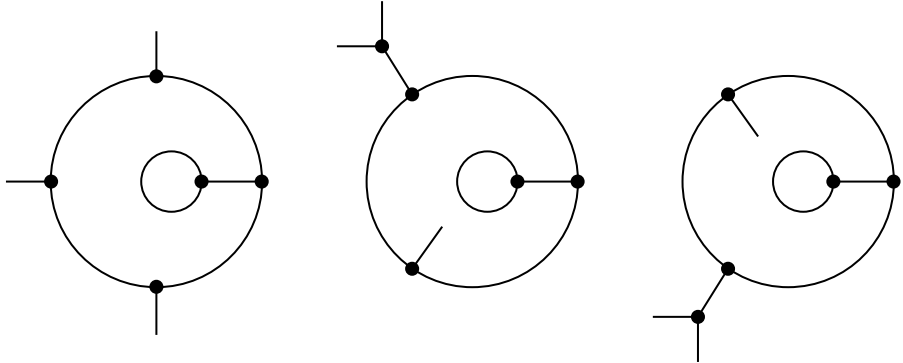}
\caption{\small Three maps with the passport $(3^5,2^61^3,7^21^1)$ and with the 
monodromy group ${\rm A}_{15}$.}
\label{fig:a_15}
\end{center}
\end{figure}

Although we have concentrated here on Hurwitz groups of the form ${\rm PSL}_2(q)$, 
there are many other examples of Hurwitz groups (see~\cite{Conder-90, Conder-10} 
for detailed surveys). The smallest alternating group which is a Hurwitz group 
is ${\rm A}_{15}$. Figure~\ref{fig:a_15} gives its $(3,2,7)$-presentations, 
using the natural representation: as for all alternating groups ${\rm A}_n$ 
with $n>6$, this is the unique non-trivial representation of least degree, 
and by the simplicity of ${\rm A}_n$ it is faithful. There are three maps 
representing ${\rm A}_{15}$. The genus of the corresponding regular Hurwitz map is
$$
g \eq \frac{|{\rm A}_{15}|}{84}+1 \eq \frac{15!}{2\cdot 84}+1 \eq 7\,783\,776\,001.
$$

More generally, Conder \cite{Conder-80} has shown that all the alternating groups
${\rm A}_n$ for $n\ge 168$ are Hurwitz groups, while for 
$n\le 167$ there are exactly 64 exceptions.

An elementary  case-by-case analysis, which is tedious both to perform and to describe, shows that there are just twelve $(3,2,7)$-maps with two faces of degree~7. These are the three with monodromy group $G={\rm PSL}_2(13)$ in Figure~9, the six in Figure~10 (the two on the right with $G={\rm PSL}_3(2)$, the other four with $|G|=1344$), and the three in Figure~16 with $G={\rm A}_{15}$. Similarly (and this is much easier to see) the only such maps with just one face of degree 7 are the two trees in Figure~2 with $G={\rm PSL}_3(2)$, and one map each in Figures~4 and 7 with $G={\rm PSL}_2(7)$ and ${\rm PSL}_2(8)$. As a by-product of this analysis we have also proved that any group which has a faithful quotient of degree $n = 16, 17, ..., 20$, including  ${\rm A}_n$, is not a Hurwitz group: one needs $n\ge 21$ in order to have three faces of degree~7.

%%%%%%%%%%%%%%%%%%%%

\section{Genus $17$ revisited}\label{17again}

In Section~6 we considered the Hurwitz group $G$ of genus $17$ and order $1344$, an extension of a normal subgroup $T\cong (C_2)^3$ by ${\rm GL}_3(2)$. We proved that this extension does not split, so that $G$ is not isomorphic to the obvious group with this normal structure, namely ${\rm AGL}_3(2)$. Since the main aim of this paper is to discuss useful methods, rather than results, we outline here an alternative way of seeing this, which may be applicable in other situations.

We use the standard result~\cite[Theorem 3.12]{Brown-82} that the equivalence classes of extensions of an abelian normal subgroup $A$ by a group $Q$, with a given action of $Q$ by conjugation on $A$, correspond to the elements of the second cohomology group $H^2(Q,A)$, with the semi-direct product corresponding to the zero element. In our case $Q = {\rm GL}_3(2) = {\rm SL}_3(2)$ and the normal subgroup $A = T$ can be regarded as its natural module $({\mathbb F}_2)^3$. Now Bell~\cite{Bell-78} has computed the cohomology of the groups ${\rm SL}_n(q)$ on their natural modules, and in this case $|H^2(Q,A)| = 2$, proving the existence of a non-split extension. We also have $Q\cong {\rm PSL}_2(7)$, and the presentation of this group in~\cite[\S7.5]{CM-80} shows that  $T$ is the normal closure in $G$ of the element $t := (yz^3)^4$, with a basis consisting of the commuting involutions $t, t^x$ and $ t^{x^2}$. We have given specific permutations in ${\rm S}_{14}$ for $x, y$ and $z$, so in principle one can compute the cocycle corresponding to this extension, and check that it is not a coboundary, proving that our extension $G$ does not split.

%%%%%%%%%%%%%%%%%%%%%%%

\ssk

\noindent School of Mathematics \hskip100pt LaBRI

\noindent University of Southampton \hskip 79pt Universit\'e de Bordeaux 

\noindent Highfield  \hskip 166pt 351 Cours de la Lib\'eration

\noindent Southampton SO17 1BJ \hskip 92pt Talence Cedex F-33405 

\noindent UK \hskip 192pt France

\msk

\noindent{\tt G.A.Jones@maths.soton.ac.uk} \hskip 56pt {\tt zvonkin@labri.fr}

\end{document}